\documentclass[11pt,twoside]{amsart}
\usepackage{graphicx}
\usepackage{amssymb}
\usepackage{amsmath}
\usepackage{amsthm}
\usepackage{amsfonts}
\usepackage{amsmath, amsfonts, amssymb}
\usepackage{multirow}
\usepackage{rotating}
\usepackage{multirow}
\usepackage{rotating}
\usepackage{multirow}
\usepackage{rotating}
\numberwithin{equation}{section}

\begin{document}
\begin{center}\small{In the name of Allah, Most gracious, Most Merciful.}\end{center}
\vspace{2cm}
\title{Complete Classification of Two-Dimensional Algebras}

\author{H.Ahmed$^1$, U.Bekbaev$^2$, I.Rakhimov$^3$}

\thanks{{\scriptsize
		emails: $^1$houida\_m7@yahoo.com; $^2$bekbaev@iium.edu.my;
		$^3$rakhimov@upm.edu.my.}}
\maketitle
\begin{center}
	{\scriptsize \address{$^{1}$Department of Math., Faculty of Science, UPM,
			Selangor, Malaysia $\&$ \\ Depart. of Math., Faculty of Science, Taiz University,
			Taiz, Yemen}\\
		\address{$^2$Department of Science in Engineering, Faculty of
			Engineering, IIUM, Malaysia}\\
		\address{$^3$Laboratory of Cryptography, Analysis and Structure,\\
			Institute for Mathematical Research (INSPEM), UPM, Serdang, Selangor, Malaysia}}
\end{center}
\begin{abstract}
 A complete classification of two-dimensional algebras over algebraically closed fields is provided.  \\
   
  \end{abstract}
  
\maketitle
\section{Introduction}

 The classification problem of finite dimensional algebras is one of the important problems of modern algebra. So far two approaches are known to the solution of the problem. One of them is structural (basis free, invariant) approach. For instance, the classification of finite dimensional simple and semi-simple associative algebras by Wedderburn and simple and semi-simple Lie algebras by Cartan are examples of such approach. But it is observed that this approach becomes more difficult when one considers more general types of algebras. Another approach to the solution of the problem is coordinate based approach (see \cite{B2015,D2003,E2004,G2011,P2011} for the latest results). These two approaches somehow are complementary to each other.
 
 In two-dimensional case a complete classification, by basis free approach, is stated in \cite{P2000} over any basic field. In this paper we follow the coordinate (basis, structure constants) based approach, we consider such problem for two-dimensional algebras over any algebraically closed field $\mathbf{F}$. We provide a list of algebras, given by their matrices of structure constants (MSC), such that any 2-dimensional algebra is isomorphic to only one algebra from the list. A similar result has been stated in \cite{G2011}. The main difference between these two results is as follows. In \cite{G2011} the authors can state the existence only whereas the uniqueness can not be guaranteed. Moreover, the approach to the problem followed in the present paper is totally different that in \cite{G2011}. Some details of our approach for arbitrary dimensional case have been given in \cite{B2015}.
 
 The paper is organized as follows. First we study the problem over algebraically closed fields of characteristics not 2 and 3, then the solution is given for algebras over algebraically closed fields of characteristics 2 and the last subsection of paper contains the result over algebraically closed fields of characteristics 3. In each of these cases we present the list of algebras via their matrices of structure constants.
 
\section{Classification of two-dimensional algebras}

To classify the main part of two-dimensional algebras we use a particular case of the following result from \cite{B2015}.  Let $n$, $m$  be any natural numbers, $\tau: (G,V)\rightarrow V$ be a fixed linear algebraic representation of an algebraic subgroup $G$ of $GL(m,\mathbf{F})$ on the $n$-dimensional vector space $V$ over $\mathbf{F}$.  Assume that there exists a nonempty $G$-invariant subset $V_0$ of $V$ and an algebraic map $P: V_0\rightarrow G$ such that
\begin{center}
	$P(\tau( g,\mathbf{v}))=P(\mathbf{v}) g^{-1}$,
\end{center}
whenever  $\mathbf{v}\in V_0$ and $g\in G$. The following result holds true \cite{B2015}.

{\bf Theorem.} \textit{ Elements $\mathbf{u},\mathbf{v}\in V_0$ are $G$-equivalent, that is $\mathbf{u}=\tau( g,\mathbf{v})$ for some $g\in G$, if and only if \[\tau(
	P(\mathbf{u}),\mathbf{u})=\tau( P(\mathbf{v}),\mathbf{v}).\]}

Let $\mathbf{A}$ be any 2 dimensional algebra over $\mathbf{F}$ with multiplication $\cdot$ given by a bilinear map $(\mathbf{u},\mathbf{v})\mapsto \mathbf{u}\cdot \mathbf{v}$ whenever $\mathbf{u}, \mathbf{v}\in \mathbf{A}$. If $e=(e^1,e^2)$ is 
basis for $\mathbf{A}$ as a vector space over $\mathbf{F}$ then one can represent this bilinear map by a matrix \[A=\left(
\begin{array}{cccc}
A^{1}_{1,1} & A^{1}_{1,2}& A^{1}_{2,1}& A^{1}_{2,2} \\
A^{2}_{1,1} & A^{2}_{1,2}& A^{2}_{2,1}& A^{2}_{2,2} \\
\end{array}
\right) \in Mat(2\times 4;\mathbf{F})\] such that \[\mathbf{u}\cdot \mathbf{v}=eA(u\otimes v)\] for any $\mathbf{u}=eu,\mathbf{v}=ev,$
where $u=(u_1, u_2),$ and  $v=(v_1, v_2)$ are column coordinate vectors of $\mathbf{u}$ and $\mathbf{v},$ respectively, $(u\otimes v)=(u_1v_1,u_1v_2,u_2v_1,u_2v_2)$, $e^i\cdot e^j=A^{1}_{i,j}e^1+A^{2}_{i,j}e^2$ whenever $i,j=1,2$.
So the algebra $\mathbf{A}$ is presented by the matrix $A\in Mat(2\times 4;\mathbf{F})$ (called the matrix of MSC of $\mathbf{A}$ with respect to the basis $e$).

If $e'=(e'^1,e'^2)$ is also another basis for $\mathbf{A}$, $g\in G=GL(2,\mathbf{F})$, $e'g=e$ and  $\mathbf{u}\cdot \mathbf{v}=e'B(u'\otimes v')$, where $\mathbf{u}=e'u',\mathbf{v}=e'v'$, then
\[\mathbf{u}\cdot \mathbf{v}=eA(u\otimes v)=e'B(u'\otimes v')=eg^{-1}B(gu\otimes gv)=eg^{-1}B(g\otimes g)(u\otimes v)\] as far as $\mathbf{u}=eu=e'u'=eg^{-1}u',\mathbf{v}=ev=e'v'=eg^{-1}v'$.
Therefore the equality
\[B=gA(g^{-1})^{\otimes 2}\] is valid, where for $g^{-1}=\left(\begin{array}{cccc} \xi_1& \eta_1\\ \xi_2& \eta_2\end{array}\right)$ one has
\[(g^{-1})^{\otimes 2}=g^{-1}\otimes g^{-1}=\left(\begin{array}{cccc} \xi _1^2 & \xi _1 \eta _1 & \xi _1 \eta _1 & \eta _1^2 \\
\xi _1 \xi _2 & \xi _1 \eta _2 & \xi _2 \eta _1 & \eta _1 \eta _2 \\
\xi _1 \xi _2 & \xi _2 \eta _1 & \xi _1 \eta _2 & \eta _1 \eta _2 \\
\xi _2^2 & \xi _2 \eta _2 & \xi _2 \eta _2 & \eta _2^2
\end{array}
\right).\]

{\bf Definition.} Two-dimensional algebras $\mathbf{A}$, $\mathbf{B}$, given by
	their matrices of structural constants $A$, $B$, are said to be isomorphic if $B=gA(g^{-1})^{\otimes 2}$ holds true for some $g\in GL(2,\mathbf{F})$.

Note that the following identities \begin{equation} \label{2.1} Tr_1(gA(g^{-1})^{\otimes 2})=Tr_1(A)g^{-1},\ \ Tr_2(gA(g^{-1})^{\otimes 2})=Tr_2(A)g^{-1}\end{equation} hold true whenever $A\in Mat(2\times 4,\mathbf{F}), g\in GL(2,\mathbf{F})$, where \[Tr_1(A)=(A^1_{1,1}+A^2_{2,1}, A^1_{1,2}+A^2_{2,2}),\ Tr_2(A)=(A^1_{1,1}+A^2_{1,2}, A^1_{2,1}+A^2_{2,2})\] are row vectors.

We divide $Mat(2\times 4,\mathbf{F})$ into the following five disjoint subsets:
\begin{itemize}
	\item[1.] All $A$ for which the system $\{Tr_1(A),Tr_2(A)\}$ is linear independent.
	\item[2.] All $A$ for which the system $\{Tr_1(A),Tr_2(A)\}$ is linear dependent and $Tr_1(A),Tr_2(A)$ are nonzero vectors.
	\item[3.] All $A$ for which $Tr_1(A)$ is nonzero vector and $Tr_2(A)=(0,0)$.
	\item[4.] All $A$ for which $Tr_1(A)=(0,0)$ and $Tr_2(A)$ is nonzero vector.
	\item[5.] All $A$ for which$Tr_1(A)=Tr_2(A)=(0,0)$.
\end{itemize}

Due to (1) it is clear that algebras with matrices from these different classes can not be isomorphic. We deal with each of these subsets separately.
Further, for the simplicity, we use the notation \[A=\left(\begin{array}{cccc} \alpha_1 & \alpha_2 & \alpha_3 &\alpha_4\\ \beta_1 & \beta_2 & \beta_3 &\beta_4\end{array}\right),\] where
$\alpha_1, \alpha_2, \alpha_3, \alpha_4, \beta_1, \beta_2, \beta_3, \beta_4$ stand for any elements of $\mathbf{F}$.

As it has been stated above we consider three cases
corresponding to $Char(\mathbf{F})\neq 2,3,$ $Char(\mathbf{F})=2$ and $Char(\mathbf{F})=3,$ respectively.

\subsection{Characteristic not $2$, $3$ case}

{\bf Theorem 1.} \textit{Over an algebraically closed field $\mathbf{F}$ with the characteristic not 2 and 3, any non-trivial 2-dimensional algebra is isomorphic to only one of the following listed, by their matrices of structure constants,  algebras:
	\[A_{1}(\mathbf{c})=\left(
	\begin{array}{cccc}
	\alpha_1 & \alpha_2 &\alpha_2+1 & \alpha_4 \\
	\beta_1 & -\alpha_1 & -\alpha_1+1 & -\alpha_2
	\end{array}\right),\ \mbox{where}\ \mathbf{c}=(\alpha_1, \alpha_2, \alpha_4, \beta_1)\in \mathbf{F}^4,\]
	\[A_{2}(\mathbf{c})=\left(
	\begin{array}{cccc}
	\alpha_1 & 0 & 0 & 1 \\
	\beta _1& \beta _2& 1-\alpha_1&0
	\end{array}\right)\simeq \left(
	\begin{array}{cccc}
	\alpha_1 & 0 & 0 & 1 \\
	-\beta _1& \beta _2& 1-\alpha_1&0
	\end{array}\right),\ \mbox{where}\ \mathbf{c}=(\alpha_1, \beta_1, \beta_2)\in \mathbf{F}^3,\]
	\[A_{3}(\mathbf{c})=\left(
	\begin{array}{cccc}
	0 & 1 & 1 & 0 \\
	\beta _1& \beta _2 & 1&-1
	\end{array}\right),\ \mbox{where}\ \mathbf{c}=(\beta_1, \beta_2)\in \mathbf{F}^2,\]
	\[A_{4}(\mathbf{c})=\left(
	\begin{array}{cccc}
	\alpha _1 & 0 & 0 & 0 \\
	0 & \beta _2& 1-\alpha _1&0
	\end{array}\right),\ \mbox{where}\ \mathbf{c}=(\alpha_1, \beta_2)\in \mathbf{F}^2,\]
	\[A_{5}(\mathbf{c})=\left(
	\begin{array}{cccc}
	\alpha_1& 0 & 0 & 0 \\
	1 & 2\alpha_1-1 & 1-\alpha_1&0
	\end{array}\right),\ \mbox{where}\ \mathbf{c}=\alpha_1\in \mathbf{F},\]
	\[A_{6}(\mathbf{c})=\left(
	\begin{array}{cccc}
	\alpha_1 & 0 & 0 & 1 \\
	\beta _1& 1-\alpha_1 & -\alpha_1&0
	\end{array}\right)\simeq \left(
	\begin{array}{cccc}
	\alpha_1 & 0 & 0 & 1 \\
	-\beta _1& 1-\alpha_1 & -\alpha_1&0
	\end{array}\right),\ \mbox{where}\ \mathbf{c}=(\alpha_1, \beta_1)\in \mathbf{F}^2,\]
	\[A_{7}(\mathbf{c})=\left(
	\begin{array}{cccc}
	0 & 1 & 1 & 0 \\
	\beta_1& 1& 0&-1
	\end{array}\right),\ \mbox{where}\ \mathbf{c}=\beta_1\in \mathbf{F},\]
	\[A_{8}(\mathbf{c})=\left(
	\begin{array}{cccc}
	\alpha_1 & 0 & 0 & 0 \\
	0 & 1-\alpha_1 & -\alpha_1&0
	\end{array}\right),\ \mbox{where}\ \mathbf{c}=\alpha_1\in \mathbf{F},\]
	\[A_{9}=\left(
	\begin{array}{cccc}
	\frac{1}{3}& 0 & 0 & 0 \\
	1 & \frac{2}{3} & -\frac{1}{3}&0
	\end{array}\right),\
	\ A_{10}=\left(
	\begin{array}{cccc}
	0 & 1 & 1 & 0 \\
	0 &0& 0 &-1
	\end{array}
	\right),\]
	\[ A_{11}=\left(
	\begin{array}{cccc}
	0 & 1 & 1 & 0 \\
	1 &0& 0 &-1
	\end{array}
	\right),\
	\ A_{12}=\left(
	\begin{array}{cccc}
	0 & 0 & 0 & 0 \\
	1 &0&0 &0\end{array}
	\right).\]}

{\bf Proof.}
{\bf The first subset case.} In this case let $P(A)$ stand for a nonsingular matrix $\left(\begin{array}{cc} \alpha_1+\beta_3& \alpha_2+\beta_4\\ \alpha_1+\beta_2 & \alpha_3+\beta_4\end{array}\right)$  with rows
\[Tr_1(A)=(\alpha _1+\beta _3, \alpha _2+\beta _4),\ Tr_2(A)=(\alpha _1+\beta _2, \alpha _3+\beta _4).\] Due to (\ref{2.1}) we have the identity \[P(g A(g^{-1})^{\otimes 2})= P(A)g^{-1},\] which means that one can apply the  Theorem as $V=\mathbf{F}^8$ and $\tau(g,A)= g A(g^{-1})^{\otimes 2}.$ In this case for $V_0$ one can take \[ V_0=\{A\in Mat(2\times 4,\mathbf{F}): \det(P(A))=\alpha_1(\alpha_3-\alpha_2)+\beta_4(\beta_3-\beta_2)+\alpha_3\beta_3-\alpha_2\beta_2\neq 0.\}\] Therefore,
two-dimensional algebras $\mathbf{A}$, $\mathbf{B}$, given by
their matrices of structure constants $A$, $B \in V_0$, are isomorphic if and only if the equality \[P(B)B(P(B)^{-1}\otimes P(B)^{-1})=P(A)A(P(A)^{-1}\otimes P(A)^{-1})\] holds true.

For $\xi_1=\frac{\alpha_3+\beta_4}{\Delta}$, $\xi_2=-\frac{\alpha_1+\beta_2}{\Delta}$, $\eta_1=-\frac{\alpha_2+\beta_4}{\Delta}$, $\eta_2=\frac{\alpha_1+\beta_3}{\Delta}$, where $\Delta=\Delta(A)=\det(P(A))$, one has
\[P(A)^{-1}=\left(\begin{array}{cc} \xi_1& \eta_1\\ \xi_2& \eta_2\end{array}\right)\] and
$A'=\left(\begin{array}{cccc} \alpha'_1 & \alpha'_2 & \alpha'_3 &\alpha'_4\\ \beta'_1 & \beta'_2 & \beta'_3 &\beta'_4\end{array}\right)=P(A)A(P(A)^{-1}\otimes P(A)^{-1})$ is matrix consisting of columns
\[\left(\begin{array}{c}
-\frac{\beta _1 \eta _1 \xi _1^2}{\Delta }+\frac{\alpha _1 \eta _2 \xi _1^2}{\Delta }+\frac{2 \alpha _1 \eta _1 \xi _1 \xi _2}{\Delta }-\frac{2 \beta _4 \eta _2 \xi _1 \xi _2}{\Delta }-\frac{2 \eta _1 \eta _2 \xi _1 \xi _2}{\Delta ^2}+\frac{\eta _2 \xi _1^2 \xi _2}{\Delta ^2}-\frac{\beta _4 \eta _1 \xi _2^2}{\Delta }+\frac{\alpha _4 \eta _2 \xi _2^2}{\Delta }+\frac{\eta _1 \xi _1 \xi _2^2}{\Delta ^2}  \\
\frac{\beta _1 \xi _1^3}{\Delta }-\frac{3 \alpha _1 \xi _1^2 \xi _2}{\Delta }+\frac{\eta _2 \xi _1^2 \xi _2}{\Delta ^2}+\frac{3 \beta _4 \xi _1 \xi _2^2}{\Delta }+\frac{\eta _1 \xi _1 \xi _2^2}{\Delta ^2}-\frac{2 \xi _1^2 \xi _2^2}{\Delta ^2}-\frac{\alpha _4 \xi _2^3}{\Delta }
\end{array}
\right),\]
\[\left(\begin{array}{c}-\frac{\beta _1 \eta _1^2 \xi _1}{\Delta }+\frac{2 \alpha _1 \eta _1 \eta _2 \xi _1}{\Delta }-\frac{\beta _4 \eta _2^2 \xi _1}{\Delta }-\frac{\eta _1 \eta _2^2 \xi _1}{\Delta ^2}+\frac{\alpha _1 \eta _1^2 \xi _2}{\Delta }-\frac{2 \beta _4 \eta _1 \eta _2 \xi _2}{\Delta }-\frac{\eta _1^2 \eta _2 \xi _2}{\Delta ^2}+\frac{\alpha _4 \eta _2^2 \xi _2}{\Delta }+\frac{2 \eta _1 \eta _2 \xi _1 \xi _2}{\Delta ^2}\\
\frac{\beta _1 \eta _1 \xi _1^2}{\Delta }-\frac{\alpha _1 \eta _2 \xi _1^2}{\Delta }-\frac{2 \alpha _1 \eta _1 \xi _1 \xi _2}{\Delta }+\frac{2 \beta _4 \eta _2 \xi _1 \xi _2}{\Delta }+\frac{2 \eta _1 \eta _2 \xi _1 \xi _2}{\Delta ^2}-\frac{\eta _2 \xi _1^2 \xi _2}{\Delta ^2}+\frac{\beta _4 \eta _1 \xi _2^2}{\Delta }-\frac{\alpha _4 \eta _2 \xi _2^2}{\Delta }-\frac{\eta _1 \xi _1 \xi _2^2}{\Delta ^2}
\end{array}
\right),\]

\[\left(\begin{array}{c} -\frac{\beta _1 \eta _1^2 \xi _1}{\Delta }+\frac{2 \alpha _1 \eta _1 \eta _2 \xi _1}{\Delta }-\frac{\beta _4 \eta _2^2 \xi _1}{\Delta }-\frac{\eta _1 \eta _2^2 \xi _1}{\Delta ^2}+\frac{\eta _2^2 \xi _1^2}{\Delta ^2}+\frac{\alpha _1 \eta _1^2 \xi _2}{\Delta }-\frac{2 \beta _4 \eta _1 \eta _2 \xi _2}{\Delta }-\frac{\eta _1^2 \eta _2 \xi _2}{\Delta ^2}+\frac{\alpha _4 \eta _2^2 \xi _2}{\Delta }+\frac{\eta _1^2 \xi _2^2}{\Delta ^2}\\
\frac{\beta _1 \eta _1 \xi _1^2}{\Delta }-\frac{\alpha _1 \eta _2 \xi _1^2}{\Delta }+\frac{\eta _2^2 \xi _1^2}{\Delta ^2}-\frac{2 \alpha _1 \eta _1 \xi _1 \xi _2}{\Delta }+\frac{2 \beta _4 \eta _2 \xi _1 \xi _2}{\Delta }-\frac{\eta _2 \xi _1^2 \xi _2}{\Delta ^2}+\frac{\beta _4 \eta _1 \xi _2^2}{\Delta }+\frac{\eta _1^2 \xi _2^2}{\Delta ^2}-\frac{\alpha _4 \eta _2 \xi _2^2}{\Delta }-\frac{\eta _1 \xi _1 \xi _2^2}{\Delta ^2}\end{array}
\right),\]

\[\left(\begin{array}{c} -\frac{\beta _1 \eta _1^3}{\Delta }+\frac{3 \alpha _1 \eta _1^2 \eta _2}{\Delta }-\frac{3 \beta _4 \eta _1 \eta _2^2}{\Delta }-\frac{2 \eta _1^2 \eta _2^2}{\Delta ^2}+\frac{\alpha _4 \eta _2^3}{\Delta }+\frac{\eta _1 \eta _2^2 \xi _1}{\Delta ^2}+\frac{\eta _1^2 \eta _2 \xi _2}{\Delta ^2}\\
\frac{\beta _1 \eta _1^2 \xi _1}{\Delta }-\frac{2 \alpha _1 \eta _1 \eta _2 \xi _1}{\Delta }+\frac{\beta _4 \eta _2^2 \xi _1}{\Delta }+\frac{\eta _1 \eta _2^2 \xi _1}{\Delta ^2}-\frac{\alpha _1 \eta _1^2 \xi _2}{\Delta }+\frac{2 \beta _4 \eta _1 \eta _2 \xi _2}{\Delta }+\frac{\eta _1^2 \eta _2 \xi _2}{\Delta ^2}-\frac{\alpha _4 \eta _2^2 \xi _2}{\Delta }-\frac{2 \eta _1 \eta _2 \xi _1 \xi _2}{\Delta ^2}\end{array}
\right).\]

From the above it can be easily seen that $\alpha'_1=-\beta'_2$, $\alpha'_3=\alpha'_2+1$, $\beta'_3=\beta'_2+1$ and $\beta'_4=-\alpha'_2,$ i.e.
\[A'=\left(
\begin{array}{cccc}
\alpha'_1 & \alpha'_2 &\alpha'_2+1 & \alpha'_4 \\
\beta'_1 & -\alpha'_1 & -\alpha'_1+1 & -\alpha'_2
\end{array}\right).\]
Therefore, the main role in the finding of $A'$ is played by the functions

$\alpha'_1=-\frac{\beta _1 \eta _1 \xi _1^2}{\Delta }+\frac{\alpha _1 \eta _2 \xi _1^2}{\Delta }+\frac{2 \alpha _1 \eta _1 \xi _1 \xi _2}{\Delta }-\frac{2 \beta _4 \eta _2 \xi _1 \xi _2}{\Delta }-\frac{2 \eta _1 \eta _2 \xi _1 \xi _2}{\Delta ^2}+\frac{\eta _2 \xi _1^2 \xi _2}{\Delta ^2}-\frac{\beta _4 \eta _1 \xi _2^2}{\Delta }+\frac{\alpha _4 \eta _2 \xi _2^2}{\Delta }+\frac{\eta _1 \xi _1 \xi _2^2}{\Delta ^2},$

$\beta'_1=\frac{\beta _1 \xi _1^3}{\Delta }-\frac{3 \alpha _1 \xi _1^2 \xi _2}{\Delta }+\frac{\eta _2 \xi _1^2 \xi _2}{\Delta ^2}+\frac{3 \beta _4 \xi _1 \xi _2^2}{\Delta }+\frac{\eta _1 \xi _1 \xi _2^2}{\Delta ^2}-\frac{2 \xi _1^2 \xi _2^2}{\Delta ^2}-\frac{\alpha _4 \xi _2^3}{\Delta }
,$

$\alpha'_2=-\frac{\beta _1 \eta _1^2 \xi _1}{\Delta }+\frac{2 \alpha _1 \eta _1 \eta _2 \xi _1}{\Delta }-\frac{\beta _4 \eta _2^2 \xi _1}{\Delta }-\frac{\eta _1 \eta _2^2 \xi _1}{\Delta ^2}+\frac{\alpha _1 \eta _1^2 \xi _2}{\Delta }-\frac{2 \beta _4 \eta _1 \eta _2 \xi _2}{\Delta }-\frac{\eta _1^2 \eta _2 \xi _2}{\Delta ^2}+\frac{\alpha _4 \eta _2^2 \xi _2}{\Delta }+\frac{2 \eta _1 \eta _2 \xi _1 \xi _2}{\Delta ^2},$

$\alpha'_4=-\frac{\beta _1 \eta _1^3}{\Delta }+\frac{3 \alpha _1 \eta _1^2 \eta _2}{\Delta }-\frac{3 \beta _4 \eta _1 \eta _2^2}{\Delta }-\frac{2 \eta _1^2 \eta _2^2}{\Delta ^2}+\frac{\alpha _4 \eta _2^3}{\Delta }+\frac{\eta _1 \eta _2^2 \xi _1}{\Delta ^2}+\frac{\eta _1^2 \eta _2 \xi _2}{\Delta ^2}.$

It can be easily verified that these functions can take any values in $\mathbf{F}.$ Therefore, the values of $\alpha'_1,\beta'_1, \alpha'_2, \alpha'_4$ also can be any elements in $\mathbf{F}$. Using obvious redenotions one can list the following non-isomorphic ``canonical'' algebras from the first subset, given by their MSCs as follows:
\[ A_{1}(\mathbf{c})=\left(
\begin{array}{cccc}
\alpha_1 & \alpha_2 &\alpha_2+1 & \alpha_4 \\
\beta_1 & -\alpha_1 & -\alpha_1+1 & -\alpha_2
\end{array}\right),\ \mbox{where}\ \mathbf{c}=(\alpha_1, \alpha_2, \alpha_4, \beta_1) \in \mathbf{F}^4.\]

For any algebra $\mathbf{A}$ from the first subset there exists a unique algebra from  the class $A_{1}(\mathbf{c})$ isomorphic to $\mathbf{A}$. Note that the appearance of $A_{1}(\mathbf{c})$ for the first subset does not depend on either algebraically closeness of $\mathbf{F}$ or its characteristic.

{\bf The second and third subset cases.} In these cases one can make $Tr_1(A)g=(1,0)$. It implies that $Tr_2(A)g=(\lambda,0)$. Therefore, we can assume that
\[Tr_1(A)=(\alpha _1+\beta _3,\alpha _2+\beta _4)=(1,0)\ \mbox{and}\ Tr_2(A)=(\lambda(\alpha _1+\beta _3),\lambda(\alpha _2+\beta _4))=(\lambda,0).\]
Here if $\lambda$ is zero we can also cover the third subset. So let us consider

\[A=\left(
\begin{array}{cccc}
\alpha _1 & \alpha _2 & \lambda  \alpha _2+(\lambda -1)\beta _4 & \alpha _4 \\
\beta _1 & (\lambda -1)\alpha _1+\lambda  \beta _3 & \beta _3 & \beta _4
\end{array}
\right)=\left(
\begin{array}{cccc}
\alpha _1 & \alpha _2 & \alpha _2 & \alpha _4 \\
\beta _1 & \lambda-\alpha _1 & 1-\alpha _1& -\alpha _2
\end{array}\right),\] with respect to $g\in GL(2,\mathbf{F})$ of the form \[g^{-1}=\left(\begin{array}{cc} 1& 0\\ \xi_2& \eta_2\end{array}\right),\] as far as
$(1,0)g^{-1}=(1,0)$ if and only if $g^{-1}$ is of the above form. In this case for the entries of $A'$ we have $A'=\left(
\begin{array}{cccc}
\alpha'_1 & \alpha'_2 & \alpha'_2 & \alpha'_4 \\
\beta'_1 & \lambda-\alpha'_1 & 1-\alpha'_1& -\alpha'_2
\end{array}\right)=gA(g^{-1})^{\otimes 2}$, where $\lambda$ is same for the $A$ and $A'$,
one has

$\alpha '_1=\frac{1}{\Delta }(\alpha _1 \eta _2 +2\alpha _2\eta _2\xi _2+\alpha _4 \eta _2 \xi _2^2)=\alpha _1+2\alpha _2\xi _2+\alpha _4\xi _2^2.$

$\alpha '_2=\frac{-1}{\Delta }(-\alpha _2\eta^2_2 -\alpha _4 \eta _2^2 \xi _2)=(\alpha _2+\alpha _4\xi _2)\eta _2.$

$\alpha '_4=\frac{-1}{\Delta }(\beta _1 \eta _1^3+[(\lambda-2)\alpha _1+(1+\lambda)\beta _3]\eta _1^2 \eta _2-[(1+\lambda)\alpha _2+(\lambda-2)\beta _4]\eta _1 \eta _2^2-\alpha _4 \eta _2^3)=\alpha _4 \eta_2^2.$

$\beta '_1=
\frac{1}{\Delta }(\beta _1 \xi _1^3+[(\lambda-2)\alpha _1 +(\lambda+1)\beta _3]\xi _1^2 \xi _2-[(1+\lambda)\alpha _2+(\lambda-2)\beta _4] \xi _1 \xi _2^2-\alpha _4 \xi _2^3)$

\qquad $=\frac{\beta _1+(1+\lambda-3 \alpha _1)\xi _2-3 \alpha _2 \xi _2^2-\alpha _4 \xi _2^3}{\eta_2}.$

So it is enough to consider the system

$\alpha '_1=\alpha _1+2 \alpha _2 \xi _2+\alpha _4 \xi _2^2,$

$\alpha '_2=(\alpha _2+\alpha_4\xi_2) \eta _2,$

$\alpha '_4=\alpha _4 \eta _2^2,$

$\beta '_1=\frac{\beta _1+(1+\lambda-3\alpha _1)\xi _2-3 \alpha _2 \xi _2^2-\alpha _4 \xi _2^3}{\eta _2}.$

Now we have to consider a few cases again.

Case 1: $\alpha_4\neq 0$. In this case one can make $$\alpha '_2=0, \alpha '_4=1, \alpha '_1=\alpha_1-\frac{\alpha^2_2}{\alpha_4}\ \mbox{and}\ \beta '_1=\sqrt{\alpha_4}\left(\beta _1-\left(1+\lambda-3\alpha _1\right)\frac{\alpha_2}{\alpha_4}-2\frac{\alpha^3_2}{\alpha^2_4}\right)$$
and
once again using redenotion one can represent the corresponding ``canonical'' MSCs as follows
\[ A_{2}(\mathbf{c})=\left(
\begin{array}{cccc}
\alpha_1 & 0 & 0 & 1 \\
\beta _1& \beta_2 & 1-\alpha_1&0
\end{array}\right),\ \mbox{where}\ \mathbf{c}=(\alpha_1, \beta_1, \beta_2)\in \mathbf{F}^3,\] as far as the expressions $\alpha_1-\frac{\alpha^2_2}{\alpha_4}$, $\sqrt{\alpha_4}(\beta _1-(1+\lambda-3\alpha _1)\frac{\alpha_2}{\alpha_4}-2\frac{\alpha^3_2}{\alpha^2_4})$, $\lambda-(\alpha_1-\frac{\alpha^2_2}{\alpha_4} ) $ may have any values in $\mathbf{F}$. It can be checked that algebras $\left(
\begin{array}{cccc}
	\alpha_1 & 0 & 0 & 1 \\
	\beta _1& \beta_2 & 1-\alpha_1&0
\end{array}\right)$ and $\left(
\begin{array}{cccc}
\alpha_1 & 0 & 0 & 1 \\
-\beta _1& \beta_2 & 1-\alpha_1&0
\end{array}\right)$ are isomorphic.
Note that the result $A_{2}(\mathbf{c})$, does not depend on the characteristic of $\mathbf{F}$.

Case 2: $\alpha _4=0.$

Subcase 2 - a): $\alpha _2\neq 0.$  If $\alpha _2\neq 0$ then one can make
$\alpha '_1=0$, $\alpha '_2=1$, $\beta '_1=\alpha_2\beta _1-\frac{2+2\lambda-3\alpha _1}{4}\alpha_1$ to get the following set of canonical matrices of structural constants
\[ A_{3}(\mathbf{c})=\left(
\begin{array}{cccc}
0 & 1 & 1 & 0 \\
\beta _1& \beta_2 & 1&-1
\end{array}\right),\ \mbox{where}\ \mathbf{c}=(\beta_1, \beta_2)\in \mathbf{F}^2.\]

Subcase 2 - b): $\alpha _2= 0.$  If $\alpha _2=0$ then $\alpha '_1=\alpha _1,$ $\alpha '_2=0,$ $\alpha '_4=0,$ $\beta '_1=\frac{\beta _1+(1+\lambda-3\alpha _1)\xi _2}{\eta _2}.$

\underline{Subsubcase: 2 - b) - 1: $1+\lambda-3\alpha _1\neq 0.$} If $1+\lambda-3\alpha _1\neq 0$, that is $\lambda-\alpha _1\neq 2\alpha_1-1$,  one can make $\beta'_1=0$ to get
\[A_{4}(\mathbf{c})=\left(
\begin{array}{cccc}
\alpha _1 & 0 & 0 & 0 \\
0 & \beta_2 & 1-\alpha _1&0
\end{array}\right),\ \mbox{where}\ \mathbf{c}=(\alpha_1, \beta_2) \in \mathbf{F}^2, \ \mbox{with}\  \beta_2\neq 2\alpha_1-1.\]

\underline{Subsubcase: 2 - b) - 2: $1+\lambda-3\alpha _1= 0.$} If $1+\lambda-3\alpha _1=0$ and $\beta_1\neq 0$ one can make $\beta'_1=1$ to get
\[A_{5}(\mathbf{c})=\left(
\begin{array}{cccc}
\alpha _1 & 0 & 0 & 0 \\
1 & 2\alpha _1-1 & 1-\alpha _1&0
\end{array}\right),\ \mbox{where}\ \mathbf{c}=\alpha_1\in F.\] If $\beta_1=0$ one has $\lambda-\alpha _1= 2\alpha_1-1$ and therefore
$A'=\left(
\begin{array}{cccc}
\alpha _1 & 0 & 0 & 0 \\
0 & 2\alpha _1-1 & 1-\alpha _1&0
\end{array}\right)$\\ which is $A_{4}(\mathbf{c})$ with $\beta_2=2\alpha_1-1$.

{\bf The fourth subset case.}
By the similar justification as in the second and the third subsets case it is enough to consider  \[ A=\left(
\begin{array}{cccc}
\alpha _1 & \alpha _2 & \alpha _2 & \alpha _4 \\
\beta _1 & 1-\alpha _1 & -\alpha _1 & -\alpha _2
\end{array}
\right)\ \mbox{and}\ \ g^{-1}=\left(\begin{array}{cc} 1& 0\\ \xi_2& \eta_2\end{array}\right),\
\mbox{where one has}\]

$\alpha'_1= \alpha _1+2 \alpha _2 \xi _2+\alpha _4 \xi _2^2 ,$

$\alpha'_2= (\alpha _2+\alpha _4\xi _2)\eta_2,$

$\alpha'_4=\alpha _4 \eta _2^2,$

$\beta'_1=\frac{\beta _1+\xi _2-3 \alpha _1 \xi _2-3 \alpha _2 \xi _2^2-\alpha _4 \xi _2^3}{\eta _2}.$

Therefore we get the canonical MSCs as follows
\[A_{6}(\mathbf{c})=\left(
\begin{array}{cccc}
\alpha_1 & 0 & 0 & 1 \\
\beta _1& 1-\alpha_1 & -\alpha_1&0
\end{array}\right)\simeq \left(
\begin{array}{cccc}
\alpha_1 & 0 & 0 & 1 \\
-\beta _1& 1-\alpha_1 & -\alpha_1&0
\end{array}\right),\ \mbox{where}\ \mathbf{c}=(\alpha_1, \beta_1) \in \mathbf{F}^2,\]
(it is $A_{2}(\mathbf{c})$, where the $2^{nd}$ and $3^{rd}$ columns are interchanged and $\alpha_1+\beta_2=0$) or
\[ A_{7}(\mathbf{c})=\left(
\begin{array}{cccc}
0 & 1 & 1 & 0 \\
\beta _1& 1& 0&-1
\end{array}\right),\ \mbox{where}\ \mathbf{c}=\beta_1\in \mathbf{F},\]
(it is $A_{3}(\mathbf{\mathbf{c}})$, where the $2^{nd}$ and $3^{rd}$ columns are interchanged and $\beta_2=0$) or
\[A_{8}(\mathbf{c})=\left(
\begin{array}{cccc}
\alpha _1 & 0 & 0 & 0 \\
0 & 1-\alpha _1 & -\alpha _1&0
\end{array}\right),\ \mbox{where}\ \mathbf{c}=\alpha_1\in \mathbf{F},\]
(it is $A_{4}(\mathbf{c})$, where the $2^{nd}$ and $3^{rd}$ columns are interchanged and $\alpha_1+\beta_2=0$) or
\[A_{9}=\left(
\begin{array}{cccc}
\frac{1}{3} & 0 & 0 & 0 \\
1 & \frac{2}{3} & -\frac{1}{3} &0
\end{array}\right)\] (it is $A_{5}(\mathbf{c})$, where the $2^{nd}$ and $3^{rd}$ columns are interchanged and $3\alpha_1-1=0$).

{\bf The fifth subset case.} In this case
\[A=\left(
\begin{array}{cccc}
\alpha _1 & \alpha _2 & \alpha _2 & \alpha _4 \\
\beta _1 & -\alpha _1 & -\alpha _1 & -\alpha _2
\end{array}
\right),\ g^{-1}=\left(\begin{array}{cc} \xi_1& \eta_1\\ \xi_2& \eta_2\end{array}\right)\mbox{and}\
A'=\left(
\begin{array}{cccc}
\alpha' _1 & \alpha' _2 & \alpha' _2 & \alpha' _4 \\
\beta' _1 & -\alpha' _1 & -\alpha' _1 & -\alpha' _2
\end{array}
\right)\]
such that:

$\alpha' _1=\frac{1}{\Delta }\left(-\beta _1 \eta _1 \xi _1^2+\alpha _1 \eta _2 \xi _1^2+2 \alpha _1 \eta _1 \xi _1 \xi _2+2 \alpha _2 \eta _2 \xi _1 \xi _2+\alpha _2 \eta _1 \xi _2^2+\alpha _4 \eta _2 \xi _2^2\right),$

$\alpha' _2=\frac{-1}{\Delta }\left(\beta _1 \eta _1^2 \xi _1-2 \alpha _1 \eta _1 \eta _2 \xi _1-\alpha _2 \eta _2^2 \xi _1-\alpha _1 \eta _1^2 \xi _2-2 \alpha _2 \eta _1 \eta _2 \xi _2-\alpha _4 \eta _2^2 \xi _2\right),$

$\alpha' _4=\frac{-1}{\Delta }\left(\beta _1 \eta _1^3-3 \alpha _1 \eta _1^2 \eta _2-3 \alpha _2 \eta _1 \eta _2^2-\alpha _4 \eta _2^3\right),$

$\beta' _1=\frac{1}{\Delta }\left(\beta _1 \xi _1^3-3 \alpha _1 \xi _1^2 \xi _2-3 \alpha _2 \xi _1 \xi _2^2-\alpha _4 \xi _2^3\right).$

If $\alpha _4\neq 0$ by making $\frac{\eta_2}{\eta_1}$ equal to any root of the polynomial $p(t)=\beta _1-3 \alpha _1t-3 \alpha _2t^2-\alpha _4t^3$ one can make $\alpha'_4=0$. Therefore, further it is assumed that $\alpha _4= 0$.

Let us consider $g$ with $\eta _1 = 0$ to have $\alpha'_4=0.$ In this case $\Delta=\xi_1\eta_2$ and \\
$\alpha' _1=\xi _1\left(\alpha _1+2\alpha _2\frac{\xi_2}{\xi_1}\right),$\\
$\alpha' _2=\alpha _2\eta_2,$\\
$\beta' _1=\frac{ \xi _1^2}{\eta_2 }\left(\beta _1-3 \alpha _1 \frac{\xi_2}{\xi_1}-3 \alpha _2(\frac{\xi_2}{\xi_1})^2\right).$\\
Case a: $\alpha_2\neq 0$. One can consider $\frac{\xi_2}{\xi_1}=\frac{-\alpha _1}{2\alpha _2}$ to get $\alpha' _1=0, \alpha' _2=1$ and $\beta' _1=\xi^2_1\frac{3\alpha^2 _1+4\alpha _2\beta _1}{4}.$ Therefore one can make $\beta' _1$ equal to $0$ or $1,$ depending on $3\alpha^2 _1+4\alpha _2\beta _1$ to have
\[A_{10}=\left(
\begin{array}{cccc}
0 & 1 & 1 & 0 \\
0 &0& 0 &-1
\end{array}
\right)\ \mbox{or}\ A_{11}=\left(
\begin{array}{cccc}
0 & 1 & 1 & 0 \\
1 &0& 0 &-1
\end{array}
\right).\]
Case b: $\alpha_2=0$. Then $\alpha'_2=\alpha'_4=0$ and $\alpha'_1=\xi_1\alpha_1, \beta' _1=\frac{\xi^2_1}{\eta_2}\left(\beta _1-3 \alpha _1 \frac{\xi_2}{\xi_1}\right).$ Therefore if $\alpha_1\neq 0$ one can make $\alpha'_1=1, \beta' _1=0$ to get $A'=\left(
\begin{array}{cccc}
1 & 0 & 0 & 0 \\
0 &-1& -1 &0
\end{array}
\right)$ , which isomorphic to $A_{10},$  if $\alpha_1= 0$ then $\alpha'_1=0$ and one can make $\beta' _1=1$ to come to
\[A_{12}=\left(
\begin{array}{cccc}
0 & 0 & 0 & 0 \\
1 &0&0 &0\end{array}
\right)\]
A routine check in each subset case shows that the corresponding algebras presented above
are not isomorphic. 

\subsection{Characteristic $2$ case}
In this case the result can be summarized in the following form.

{\bf Theorem 2.}\textit{ Over an algebraically closed field $\mathbf{F}$ characteristic 2 any non-trivial 2-dimensional algebra is isomorphic to only one of the following listed, by their matrices of structure constants,  algebras:
\[A_{1,2}(\mathbf{c})=\left(
\begin{array}{cccc}
\alpha_1 & \alpha_2 &\alpha_2+1 & \alpha_4 \\
\beta_1 & -\alpha_1 & -\alpha_1+1 & -\alpha_2
\end{array}\right),\ \mbox{where}\ \mathbf{c}=(\alpha_1, \alpha_2, \alpha_4, \beta_1)\in \mathbf{F}^4,\]
\[A_{2,2}(\mathbf{c})=\left(
\begin{array}{cccc}
\alpha_1 & 0 & 0 & 1 \\
\beta _1& \beta_2 & 1-\alpha_1&0
\end{array}\right),\ \mbox{where}\ \mathbf{c}=(\alpha_1, \beta_1, \beta_2)\in \mathbf{F}^3,\]
\[A_{3,2}(\mathbf{c})=\left(
\begin{array}{cccc}
\alpha_1 & 1 & 1 & 0 \\
0& \beta_2 & 1-\alpha_1&1
\end{array}\right),\ \mbox{where}\ \mathbf{c}=(\alpha_1, \beta_2)\in \mathbf{F}^2,\]
\[A_{4,2}(\mathbf{c})=\left(
\begin{array}{cccc}
\alpha _1 & 0 & 0 & 0 \\
0 & \beta_2 & 1-\alpha _1&0
\end{array}\right),\ \mbox{where}\ \mathbf{c}=(\alpha_1,\beta_2)\in \mathbf{F}^2,\]
\[A_{5,2}(\mathbf{c})=\left(
\begin{array}{cccc}
\alpha_1 & 0 & 0 & 0 \\
1 & 1 & 1-\alpha_1&0
\end{array}\right),\ \mbox{where}\ \mathbf{c}=\alpha_1\in \mathbf{F},\]
\[A_{6,2}(\mathbf{c})=\left(
\begin{array}{cccc}
\alpha_1 & 0 & 0 & 1 \\
\beta _1& 1-\alpha_1 & -\alpha_1&0
\end{array}\right),\ \mbox{where}\ \mathbf{c}=(\alpha_1, \beta_1)\in \mathbf{F}^2,\]
\[A_{7,2}(\mathbf{c})=\left(
\begin{array}{cccc}
\alpha_1 & 1 & 1 & 0 \\
0& 1-\alpha_1& -\alpha_1&-1
\end{array}\right),\ \mbox{where}\ \mathbf{c}=\alpha_1\in \mathbf{F},\]
\[A_{8,2}(\mathbf{c})=\left(
\begin{array}{cccc}
\alpha_1 & 0 & 0 & 0 \\
0 & 1-\alpha_1 & -\alpha_1&0
\end{array}\right),\ \mbox{where}\ \mathbf{c}=\alpha_1\in \mathbf{F},\]
\[A_{9,2}=\left(
\begin{array}{cccc}
1 & 0 & 0 & 0 \\
1 & 0 & 1&0
\end{array}\right),\ \ A_{10,2}=\left(
\begin{array}{cccc}
0 & 1 & 1 & 0 \\
0 &0& 0 &-1
\end{array}
\right),\]
\[ A_{11,2}=\left(
\begin{array}{cccc}
1 & 1 & 1 & 0 \\
0 &-1& -1 &-1
\end{array}
\right),\
\ A_{12,2}=\left(
\begin{array}{cccc}
0 & 0 & 0 & 0 \\
1 &0&0 &0\end{array}
\right).\]}

{\bf Proof.}
In this case $\alpha=-\alpha$ for any $\alpha\in \mathbf{F}$, but sometimes we use $-\alpha$ also to keep similarity with the previous case.

{\bf The first subset case.} It has been noted earlier that in this case the result doesn't depend either on algebraically closeness of $\mathbf{F}$ or its characteristics. Therefore, one has \[A_{1,2}(\mathbf{c})=\left(
\begin{array}{cccc}
\alpha_1 & \alpha_2 &\alpha_2+1 & \alpha_4 \\
\beta_1 & -\alpha_1 & -\alpha_1+1 & -\alpha_2
\end{array}\right),\ \mbox{where}\ \mathbf{c}=(\alpha_1, \alpha_2, \alpha_4, \beta_1)\in \mathbf{F}^4.\]

{\bf The second and third subset cases.}
Obviously in this case one can make $Tr_1(A)g=(1,0)$. This implies $Tr_2(A)g=(\lambda,0)$. Therefore, it can be assumed that  $Tr_1(A)=(\alpha _1+\beta _3,\alpha _2+\beta _4)=(1,0)$ and $Tr_2(A)=( \lambda(\alpha _1+\beta _3),\lambda(\alpha _2+\beta _4)=(\lambda,0)$. Here also we allow $\lambda$ to be zero to include the third subset's case. Therefore, we consider
\[A=\left(
\begin{array}{cccc}
\alpha _1 & \alpha _2 & \alpha_2 & \alpha _4 \\
\beta _1 & \lambda-\alpha _1 & 1-\alpha _1& -\alpha _2
\end{array}\right)\ \mbox{with respect to}\ g^{-1}=\left(\begin{array}{cc} 1& 0\\ \xi_2& \eta_2\end{array}\right) \ \mbox{and}\]

$\alpha '_1=\alpha _1+\alpha _4\xi _2^2,$

$\alpha '_2=(\alpha _2+\alpha _4\xi _2)\eta _2,$

$\alpha '_4=\alpha _4 \eta_2^2,$

$\beta '_1=\frac{\beta _1+(1+\lambda-3 \alpha _1)\xi _2-3 \alpha _2 \xi _2^2-\alpha _4 \xi _2^3}{\eta_2}=\frac{\beta _1+(1+\lambda+\alpha _1)\xi _2+ \alpha _2 \xi _2^2+\alpha _4 \xi _2^3}{\eta_2}.$
\begin{itemize}
	\item[1.] As we have observed earlier in $\alpha_4\neq 0$ case one gets
	\[A_{2,2}(\mathbf{c})=\left(
	\begin{array}{cccc}
	\alpha_1 & 0 & 0 & 1 \\
	\beta _1& \beta_2 & 1-\alpha_1&0
	\end{array}\right),\ \mbox{where}\ \mathbf{c}=(\alpha_1, \beta_1, \beta_2)\in \mathbf{F}^3.\]
	\item[2.] Let now assume that $\alpha _4=0.$
	\begin{itemize}
		\item[2-a).] If $\alpha _2\neq0$ then one can make
		$\alpha '_1=\alpha_1$, $\alpha '_2=1$, $\alpha'_4=0$ and $\beta'_1=0$ to get \[A_{3,2}(\mathbf{c})=\left(
		\begin{array}{cccc}
		\alpha_1 & 1 & 1 & 0 \\
		0& \beta_2 & 1-\alpha_1&-1
		\end{array}\right), \ \ \mbox{where}\ \mathbf{c}=(\alpha_1, \beta_2)\in \mathbf{F}^2.\]
		
		\item[2-b).] If $\alpha _2=0$ then $\alpha '_1=\alpha _1,$ $\alpha '_2=0,$ $\alpha '_4=0,$ and $\beta '_1=\frac{\beta _1+(1+\lambda+\alpha _1)\xi _2}{\eta _2}.$
		\begin{itemize}
			\item[2-b)-1.] If $1+\lambda+\alpha _1\neq 0$, that is $\lambda-\alpha _1\neq 1$, one can make $\beta'_1=0$ to get
			\[A_{4,2}(\mathbf{c})=\left(
			\begin{array}{cccc}
			\alpha _1 & 0 & 0 & 0 \\
			0 & \beta_2 & 1-\alpha _1&0
			\end{array}\right),\  \ \mbox{where}\ \mathbf{c}=(\alpha_1, \beta_2)\in \mathbf{F}\ \mbox{with}\ \beta_2\neq 1.\]
			
			\item[2-b)-2.] Let $1+\lambda+\alpha _1=0$. In this case
			\begin{itemize}
				\item if $\beta_1\neq 0$ one can make $\beta'_1=1$ to have
				\[A_{5,2}(\mathbf{c})=\left(
				\begin{array}{cccc}
				\alpha _1 & 0 & 0 & 0 \\
				1 & 1 & 1-\alpha _1&0
				\end{array}\right), \ \mbox{where}\ \mathbf{c}= \alpha _1 \in \mathbf{F}.\]
				\item if $\beta_1=0$ then
				\[A'=\left(
				\begin{array}{cccc}
				\alpha _1 & 0 & 0 & 0 \\
				0 & 1 & 1-\alpha _1&0
				\end{array}\right), \ \mbox{which is}\ A_{4,2}(\mathbf{c}) \ \mbox{with}\  \beta_2=1.\]
			\end{itemize}
		\end{itemize}
	\end{itemize}
\end{itemize}
{\bf The fourth subset case.} For this case we have \[ A=\left(
\begin{array}{cccc}
\alpha _1 & \alpha _2 & \alpha _2 & \alpha _4 \\
\beta _1 & 1-\alpha _1 & -\alpha _1 & -\alpha _2
\end{array}
\right)\ \mbox{and}\ g^{-1}=\left(\begin{array}{cc} 1& 0\\ \xi_2& \eta_2\end{array}\right)\]
therefore

$\alpha'_1= \alpha _1+\alpha _4 \xi _2^2 ,$

$\alpha'_2= (\alpha _2+\alpha _4\xi _2)\eta_2,$

$\alpha'_4=\alpha _4 \eta _2^2,$

$\beta'_1=\frac{\beta _1+(1+\alpha _1)\xi _2+\alpha _2 \xi _2^2+\alpha _4 \xi _2^3}{\eta _2}.$

Hence, one gets
$A_{6,2}(\mathbf{c})=\left(
\begin{array}{cccc}
\alpha_1 & 0 & 0 & 1 \\
\beta _1& 1-\alpha_1 & -\alpha_1&0
\end{array}\right) \ \ \mbox{where}\ \mathbf{c}=(\alpha_1, \beta_1)\in \mathbf{F}^2$\\ (it is $A_{2,2}(\mathbf{c})$, where the $2^{nd}$ and $3^{rd}$ columns are interchanged and $\alpha_1+\beta_2=0$),
\[A_{7,2}(\mathbf{c})=\left(
\begin{array}{cccc}
\alpha_1 & 1 & 1 & 0 \\
0&1-\alpha_1&-\alpha_1&-1
\end{array}\right),\ \ \mbox{where}\ \mathbf{c}=\alpha_1\in \mathbf{F}\] (it is $A_{3,2}(\mathbf{c})$, where the $2^{nd}$ and $3^{rd}$ columns are interchanged and $\alpha_1+\beta_2=0$),
\[A_{8,2}(\mathbf{c})=\left(
\begin{array}{cccc}
\alpha_1 & 0 & 0 & 0 \\
0 & 1-\alpha_1 & -\alpha_1&0
\end{array}\right)\ \ \mbox{where}\ \mathbf{c}=\alpha_1\in \mathbf{F}\] (it is $A_{4,2}(\mathbf{c})$, where the $2^{nd}$ and $3^{rd}$ columns are interchanged and $\alpha_1+\beta_2=0$)
and
\[A_{9,2}=\left(
\begin{array}{cccc}
1 & 0 & 0 & 0 \\
1 &0 & 1&0
\end{array}\right)\] (it is $A_{5,2}(\mathbf{c})$, where the $2^{nd}$ and $3^{rd}$ columns are interchanged and $\alpha_1+1=0$).

{\bf The fifth subset case.} In this case
\[A=\left(
\begin{array}{cccc}
\alpha _1 & \alpha _2 & \alpha _2 & \alpha _4 \\
\beta _1 & -\alpha _1 & -\alpha _1 & -\alpha _2
\end{array}\right),\ \ g^{-1}=\left(\begin{array}{cc} \xi_1& \eta_1\\ \xi_2& \eta_2\end{array}\right)\]
and for the entries of $A'=\left(
\begin{array}{cccc}
\alpha' _1 & \alpha' _2 & \alpha' _2 & \alpha' _4 \\
\beta' _1 & -\alpha' _1 & -\alpha' _1 & -\alpha' _2
\end{array}
\right)$ one has

$\alpha' _1=\frac{1}{\Delta }\left(\beta _1 \eta _1 \xi _1^2+\alpha _1 \eta _2 \xi _1^2+ \alpha_2 \eta _1\xi^2_2+\alpha _4 \eta _2 \xi _2^2\right),$

$\alpha' _2=\frac{1}{\Delta }\left(\beta _1 \eta _1^2 \xi _1+\alpha _1 \eta _1^2 \xi _2+\alpha _2 \eta _2^2 \xi _1+\alpha _4 \eta _2^2 \xi _2\right),$

$\alpha' _4=\frac{1}{\Delta }\left(\beta _1 \eta _1^3+ \alpha _1 \eta _1^2 \eta _2+\alpha _2 \eta _1 \eta _2^2+\alpha _4 \eta _2^3\right),$

$\beta' _1=\frac{1}{\Delta }\left(\beta _1 \xi _1^3+\alpha _1 \xi _1^2 \xi _2+\alpha _2 \xi _1 \xi _2^2+\alpha _4 \xi _2^3\right).$

As it was observed one can assume that $\alpha_4=0$.

Let us consider $g$ with $\eta _1 = 0$ to have $\alpha'_4=0.$ In this case $\Delta=\xi_1\eta_2$ and \\
$\alpha' _1=\xi _1\alpha _1,$\\
$\alpha' _2=\alpha _2\eta_2,$\\
$\beta' _1=\frac{ \xi _1^2}{\eta_2 }\left(\beta _1+ \alpha _1 \frac{\xi_2}{\xi_1}+ \alpha _2(\frac{\xi_2}{\xi_1})^2\right).$\\
If one of $\alpha _1, \alpha _2$ is not zero one can make $\beta' _1=0$ and depending on $\alpha _1, \alpha _2$ one can have $\alpha' _1=0, \alpha' _2=1$ or $\alpha' _1=1, \alpha' _2=1$ or $\alpha' _1=1, \alpha' _2=0$ to get
$A_{10,2}=\left(
\begin{array}{cccc}
0 & 1 & 1 & 0 \\
0 &0& 0 &-1
\end{array}
\right)$ or $A_{11,2}=\left(
\begin{array}{cccc}
1 & 1 & 1 & 0 \\
0 &-1& -1 &-1
\end{array}
\right)$ or $A'=\left(
\begin{array}{cccc}
1 & 0 & 0 & 0 \\
0 &-1& -1 &0
\end{array}
\right)$ which is isomorphic to $A_{10,2}$\\
If $\alpha _1=\alpha _2=0$ then one can make $\beta' _1=1$ to get
\[A_{12,2}=\left(
\begin{array}{cccc}
0 & 0 & 0 & 0 \\
1 &0& 0 &0
\end{array}
\right).\]

\subsection{Characteristic $3$ case}
In this case we summarize the result as follows.

{\bf Theorem 3.}\textit{ Over an algebraically closed field $\mathbf{F}$ characteristics 3 any non-trivial 2-dimensional algebra is isomorphic to only one of the following listed, by their matrices of structure constant matrices,  algebras:
\[A_{1,3}(\mathbf{c})=\left(
\begin{array}{cccc}
\alpha_1 & \alpha_2 &\alpha_2+1 & \alpha_4 \\
\beta_1 & -\alpha_1 & -\alpha_1+1 & -\alpha_2
\end{array}\right),\ \mbox{where}\ \mathbf{c}=(\alpha_1, \alpha_2, \alpha_4, \beta_1)\in \mathbf{F}^4,\]
\[A_{2,3}(\mathbf{c})=\left(
\begin{array}{cccc}
\alpha_1 & 0 & 0 & 1 \\
\beta _1& \beta _2& 1-\alpha_1&0
\end{array}\right)\simeq\left(
\begin{array}{cccc}
\alpha_1 & 0 & 0 & 1 \\
-\beta _1& \beta _2& 1-\alpha_1&0
\end{array}\right),\ \mbox{where}\ \mathbf{c}=(\alpha_1, \beta_1, \beta_2)\in \mathbf{F}^3,\]
\[A_{3,3}(\mathbf{c})=\left(
\begin{array}{cccc}
0 & 1 & 1 & 0 \\
\beta _1& \beta _2 & 1&-1
\end{array}\right),\ \mbox{where}\ \mathbf{c}=(\beta_1, \beta_2)\in \mathbf{F}^2,\]
\[A_{4,3}(\mathbf{c})=\left(
\begin{array}{cccc}
\alpha _1 & 0 & 0 & 0 \\
0 & \beta _2& 1-\alpha _1&0
\end{array}\right),\ \mbox{where}\ \mathbf{c}=(\alpha_1, \beta_2)\in \mathbf{F}^2,\]
\[A_{5,3}(\mathbf{c})=\left(
\begin{array}{cccc}
\alpha_1& 0 & 0 & 0 \\
1 & -1-\alpha_1 & 1-\alpha_1&0
\end{array}\right),\ \mbox{where}\ \mathbf{c}=\alpha_1\in \mathbf{F},\]
\[A_{6,3}(\mathbf{c})=\left(
\begin{array}{cccc}
\alpha_1 & 0 & 0 & 1 \\
\beta _1& 1-\alpha_1 & -\alpha_1&0
\end{array}\right)\simeq \left(
\begin{array}{cccc}
\alpha_1 & 0 & 0 & 1 \\
-\beta _1& 1-\alpha_1 & -\alpha_1&0
\end{array}\right),\ \mbox{where}\ \mathbf{c}=(\alpha_1, \beta_1)\in \mathbf{F}^2,\]
\[A_{7,3}(\mathbf{c})=\left(
\begin{array}{cccc}
0 & 1 & 1 & 0 \\
\beta_1& 1& 0&-1
\end{array}\right),\ \mbox{where}\ \mathbf{c}=\beta_1\in \mathbf{F},\]
\[A_{8,3}(\mathbf{c})=\left(
\begin{array}{cccc}
\alpha_1 & 0 & 0 & 0 \\
0 & 1-\alpha_1 & -\alpha_1&0
\end{array}\right),\ \mbox{where}\ \mathbf{c}=\alpha_1\in \mathbf{F},\]
\[A_{9,3}=\left(
\begin{array}{cccc}
0 & 1& 1& 0 \\
1 &0&0 &-1\end{array}
\right),\
\ A_{10,3}=\left(
\begin{array}{cccc}
0 & 1 & 1 & 0 \\
0 &0&0 &-1\end{array}
\right),\]
\[ A_{11,3}=\left(
\begin{array}{cccc}
1 & 0 & 0 & 0 \\
1 &-1&-1 &0\end{array}
\right),\ \ A_{12,3}=\left(
\begin{array}{cccc}
0 &0 &0 & 0 \\
1 &0& 0 &0
\end{array}
\right).\]}

{\bf Proof.}
\quad {\bf The first subset case.} One has
$$A_{1,3}(\mathbf{c})=\left(
\begin{array}{cccc}
\alpha_1 & \alpha_2 &\alpha_2+1 & \alpha_4 \\
\beta_1 & -\alpha_1 & -\alpha_1+1 & -\alpha_2
\end{array}\right), \ \mbox{where}\ \mathbf{c}=(\alpha_1, \alpha_2, \alpha_4, \beta_1)\in \mathbf{F}^4.$$

{\bf The second and third subset cases.} One can consider
\[A=\left(
\begin{array}{cccc}
\alpha _1 & \alpha _2 & \alpha _2 & \alpha _4 \\
\beta _1 & \lambda-\alpha _1 & 1-\alpha _1& -\alpha _2
\end{array}\right),\] where $\lambda$ may be zero as well, with respect to
$g^{-1}=\left(\begin{array}{cc} 1& 0\\ \xi_2& \eta_2\end{array}\right)$ in that case for the entries of $A'$ we get

$\alpha '_1=\alpha _1+2 \alpha _2 \xi _2+\alpha _4 \xi _2^2,$

$\alpha '_2=(\alpha _2+\alpha_4\xi_2) \eta _2,$

$\alpha '_4=\alpha _4 \eta _2^2,$

$\beta '_1=\frac{\beta _1+(1+\lambda)\xi _2-\alpha _4 \xi _2^3}{\eta _2}.$
\begin{itemize}
	\item[1.] $\alpha_4\neq 0$. In this case we obtain
	$$ A_{2,3}(\mathbf{c})=\left(
	\begin{array}{cccc}
	\alpha_1 & 0 & 0 & 1 \\
	\beta _1& \beta_2 & 1-\alpha_1&0
	\end{array}\right)\simeq \left(
	\begin{array}{cccc}
	\alpha_1 & 0 & 0 & 1 \\
	-\beta _1& \beta_2 & 1-\alpha_1&0
	\end{array}\right),\ \mbox{where}\ \mathbf{c}=(\alpha_1, \beta_1, \beta_2)\in \mathbf{F}^3.$$
	
	\item[2.] $\alpha _4=0.$
	\begin{itemize}
		\item[2-a)] If $\alpha _2\neq0$ one can make $\alpha '_1=0$, $\alpha'_2=1$ to get \[ A_{3,3}(\mathbf{c})=\left(
		\begin{array}{cccc}
		0 & 1 & 1 & 0 \\
		\beta _1& \beta_2 & 1&-1
		\end{array}\right)\ \mbox{with}\ \mathbf{c}=(\beta_1,\beta_2)\in \mathbf{F}^2.\]
		
		\item[2-b)] If $\alpha _2=0$ then $\alpha '_1=\alpha _1,$ $\alpha '_2=0,$ $\alpha '_4=0,$ and $\beta '_1=\frac{\beta _1+(1+\lambda)\xi _2}{\eta _2}.$
		\begin{itemize}
			\item[2-b)-1.] Let $1+\lambda\neq 0$, that is $\lambda\neq -1$, one can make $\beta'_1=0$ to get
			\[A_{4,3}(\mathbf{c})=\left(
			\begin{array}{cccc}
			\alpha _1 & 0 & 0 & 0 \\
			0 & \beta_2 & 1-\alpha _1&0
			\end{array}\right), \ \mbox{where}\ \mathbf{c}=(\alpha_1,\beta_2)\in \mathbf{F}^2 \ \mbox{with}\  \beta_2\neq -1-\alpha_1.\]
			
			\item[2-b)-2.] Let $1+\lambda=0.$
			\begin{itemize}
				\item If $\beta_1\neq 0$ one can make $\beta'_1=1$ to get
				\[A_{5,3}(\mathbf{c})=\left(
				\begin{array}{cccc}
				\alpha_1 & 0 & 0 & 0 \\
				1 & -1-\alpha_1 & 1-\alpha_1&0
				\end{array}\right),\ \mbox{where}\ \mathbf{c}=\alpha_1\in \mathbf{F}.\]
				\item If $\beta_1=0$ one has $\lambda=-1$ and therefore
				$A'=\left(
				\begin{array}{cccc}
				\alpha _1 & 0 & 0 & 0 \\
				0 & -1-\alpha _1 & 1-\alpha _1&0
				\end{array}\right)$ which is $A_{4,3}(\mathbf{c})$ with $\beta_2=-1-\alpha_1$.
			\end{itemize}
		\end{itemize}
	\end{itemize}
\end{itemize}

{\bf The fourth subset case.} It is easy to see that in this case the result can be derived from the second and third subsets case. So we get
\[A_{6,3}(\mathbf{c})=\left(
\begin{array}{cccc}
\alpha_1 & 0 & 0 & 1 \\
\beta _1& 1-\alpha_1 & -\alpha_1&0
\end{array}\right)\simeq \left(
\begin{array}{cccc}
\alpha_1 & 0 & 0 & 1 \\
-\beta _1& 1-\alpha_1 & -\alpha_1&0
\end{array}\right),\ \mbox{where}\ \mathbf{c}=(\alpha_1,\beta_1)\in \mathbf{F}^2\]\\ (it is $A_{2,3}(\mathbf{c})$, where the $2^{nd}$ and $3^{rd}$ columns are interchanged and $\alpha_1+\beta_2=0$),
\[ A_{7,3}(\mathbf{c})=\left(
\begin{array}{cccc}
0 & 1 & 1 & 0 \\
\beta_1& 1& 0&-1
\end{array}\right), \ \mbox{where}\ \mathbf{c}=\beta_1\in \mathbf{F}\]\\ (it is $A_{3,3}(\mathbf{c})$, where the $2^{nd}$ and $3^{rd}$ columns are interchanged and $\beta_2=0$),
\[A_{8,3}(\mathbf{c})=\left(
\begin{array}{cccc}
\alpha_1 & 0 & 0 & 0 \\
0 & 1-\alpha_1 & -\alpha_1&0
\end{array}\right), \mbox{where}\ \mathbf{c}=\alpha_1\in \mathbf{F}\]\\ (it is $A_{4,3}(\mathbf{c})$, where the $2^{nd}$ and $3^{rd}$ columns are interchanged and $\alpha_1+\beta_2=0$).

{\bf The fifth subset case.}

Here it is enough to consider $\alpha_4=0$ case. Therefore we have

$\alpha' _1=\frac{1}{\Delta }\left(-\beta _1 \eta _1 \xi _1^2+\alpha _1 \eta _2 \xi _1^2-\alpha _1 \eta _1 \xi _1 \xi _2-\alpha _2 \eta _2 \xi _1 \xi _2+\alpha _2 \eta _1 \xi _2^2\right),$

$\alpha' _2=\frac{-1}{\Delta }\left(\beta _1 \eta _1^2 \xi _1+\alpha _1 \eta _1 \eta _2 \xi _1-\alpha _2 \eta _2^2 \xi _1-\alpha _1 \eta _1^2 \xi _2+\alpha _2 \eta _1 \eta _2 \xi _2\right),$

$\alpha' _4=\frac{-1}{\Delta }\beta _1 \eta _1^3,$

$\beta' _1=\frac{1}{\Delta }\beta _1 \xi _1^3.$

Making $\eta_1=0$ results in $\alpha'_1=\alpha_1 \xi_1-\alpha _2\xi_2,$ $\alpha'_4=0$, $\alpha'_2=\alpha _2\eta_2$,  $\beta'_1=\frac{\xi^2_1}{\eta_2}\beta _1.$
\begin{itemize}
	\item[a)] If $\alpha_2\neq 0$ one can make $\alpha'_1= 0$ and $\alpha'_2=1$.
	\begin{itemize}
		\item[a)-1.] If $\beta_1\neq 0$ one can reduce $\beta'_1=1$ to get
		$A_{9,3}=\left(
		\begin{array}{cccc}
		0 & 1 & 1 & 0 \\
		1 &0&0 &-1\end{array}
		\right).$
		
		\item[a)-2.] If $\beta_1= 0$ then $\beta'_1=0$ and one gets
		$A_{10,3}=\left(
		\begin{array}{cccc}
		0 &1 & 1 & 0 \\
		0 &0& 0 &-1
		\end{array}
		\right).$
	\end{itemize}
	\item[b)] If $\alpha_2= 0$ then $\alpha'_2= 0$.
	\begin{itemize}
		\item[b)-1.] If $\beta_1\neq 0$  we reduce $\beta'_1=1$ and $\alpha'_1$ may take one or zero depending on $\alpha_1$. Therefore one has
		\[A_{11,3}=\left(
		\begin{array}{cccc}
		1 &0 & 0 & 0 \\
		1 &-1&-1 &0
		\end{array}
		\right)\ \mbox{or}\ A_{12,3}=\left(
		\begin{array}{cccc}
		0 &0 &0 & 0 \\
		1 &0& 0 &0
		\end{array}
		\right).\]
		
		\item[b)-2.] If $\beta_1=0$  then $\beta'_1=0$, $\alpha'_4=0$ and one can make , $\alpha'_1=1$ to get
		\[A'=\left(
		\begin{array}{cccc}
		1 &0& 0 & 0 \\
		0 &-1& -1 &0
		\end{array}
		\right)\ \mbox{isomorphic to}\ A_{10,3}.\]
	\end{itemize}
\end{itemize}

{\bf Remark 1.} From the proof of Theorems it is clear that the same results remain be true if one assumes the existence of a root in $\mathbf{F}$ of every second and third order polynomial over $\mathbf{F}$ (instead of assumption that $\mathbf{F}$ to be algebraically closed). The results can be used for getting complete classification of different classes, for example,  associative, Jordan, alternative and etc., 2-dimensional algebras. 

{\bf Remark 2.} In \cite{A} the class $A_{3,2}(\mathbf{c})$ should be understood as it is in this paper.
 
{\bf Acknowledgment.}
  The first author thanks M.A.A. Ahmed for the help in computer computations, the second author's research is supported by FRGS14-153-0394, MOHE and the third author acknowledges MOHE for a support by grant 01-02-14-1591FR.

\vspace{0,2cm}

\bibliographystyle{aipproc}

\end{document}